\documentclass[preprint,authoryear]{elsarticle}

\biboptions{sort&compress}

\usepackage{graphicx}      
\usepackage{natbib}        
\usepackage{amssymb}
\usepackage{amsmath}
\usepackage{enumerate}
\usepackage{placeins}
\usepackage{epstopdf}
\usepackage[none]{hyphenat}
\newtheorem{thm}{Theorem}
\newtheorem{ass}{Assumption}
\newtheorem{lemma}{Lemma}
\newtheorem{define}{Definition}

\newproof{pf}{Proof}

\begin{document}
\bibliographystyle{model5-names}
\pagestyle{plain}
\setcounter{page}{0}
\pagenumbering{arabic}
\begin{frontmatter}

\title{Constraint Back-Offs for Safe, Sufficient Excitation: A General Theory with Application to Experimental Optimization} 

\author{Gene A. Bunin}
\ead{gene.bunin@ronininstitute.org}

\address{Ronin Institute for Independent Scholarship}

\begin{abstract}                

In many experimental settings, one is tasked with obtaining information about certain relationships by applying perturbations to a set of independent variables and noting the changes in the set of dependent ones. While traditional design-of-experiments methods are often well-suited for this, the task becomes significantly more difficult in the presence of constraints, which may make it impossible to sufficiently excite the experimental system without incurring constraint violations. The key contribution of this paper consists in deriving constraint back-off sizes sufficient to guarantee that one can always perturb in a ball of radius $\delta_e$ without leaving the constrained space, with $\delta_e$ set by the user. Additionally, this result is exploited in the context of experimental optimization to propose a constrained version of G.~E.~P. Box's evolutionary operation technique. The proposed algorithm is applied to three case studies and is shown to consistently converge to the neighborhood of the optimum without violating constraints.

\;

\noindent Keywords: design measures for robustness, evolutionary operation, optimization under uncertainties, experiment design

\end{abstract}

\end{frontmatter}

\section{Introduction}
\label{sec:intro}

In most branches of science, one often encounters systems where the relationship between some experimental response and a finite number of independent variables needs to be studied \citep{Myers2009, Montgomery2012}, and it is generally assumed that the response is a dependent variable and a function of the independent ones. Mathematically, the experimental quantity may be stated as the function $f : \mathbb{R}^{n_u} \rightarrow \mathbb{R}$, while ${\bf u} \in \mathbb{R}^{n_u}$ may be used to denote the vector of independent variables ${\bf u} = (u_1, ..., u_{n_u})$. One is then left with the task of identifying $f({\bf u})$. The identification may be either local or global, and is usually done by conducting a series of experiments with different values of ${\bf u}$, observing the resulting $f({\bf u})$ values, and performing some sort of regression. Such procedures are typically used to:

\begin{enumerate}[(i)]
\item construct data-driven model approximations of $f$ for when $f$ is difficult to model via first principles \citep{Jones1998, Myers2009, Montgomery2012},
\item estimate the uncertain parameters of an already available model \citep{Box1990, Chen1987, Pfaff2006, Quelhas:12},
\item explore how the function value changes so as to find conditions for which the value is minimized, maximized, or equal to a certain quantity \citep{Robbins1951, Box:51, Lewis2000, Conn2009}. 
\end{enumerate}

It is the case for many problems that the experimental space of interest is a box defined by the constraints $u_i^L \leq u_i \leq u_i^U, \; i = 1,...,n_u$, where ${\bf u}^L = \left( u_1^L, ..., u_{n_u}^L \right)$ and ${\bf u}^U = \left( u_1^U, ..., u_{n_u}^U \right)$ are the lower and upper limits on the independent variables, respectively. Such problems typically correspond to simple set-ups that do not possess major safety limitations, and where testing any variable combination in the experimental space is permissible. Obtaining knowledge about $f$ is not difficult in such conditions, and the traditional design-of-experiments techniques \citep{Montgomery2012} are perfectly appropriate here.

However, there still exists a fair share of problems -- many of them corresponding to continuous or batch chemical processes \citep{ExpOpt} -- where additional constraints enter to reduce the experimental space in a nontrivial manner. These constraints may be expressed as the $n_g$ inequalities

$$
g_{j} ({\bf u}) \leq 0, \; j = 1,...,n_g.
$$

\noindent In some problems, the functions $g_j$ may represent \emph{experimental} relationships that, like $f$, can only be divined empirically. It often happens that such experimental constraints are safety or economic limitations -- they could, for example, represent an upper limit on the temperature in a continuous reactor, or a lower limit on the purity of a batch-produced chemical. Despite the violation of such constraints being highly undesirable, or even dangerous, there currently exists no easy-to-implement, theoretically rigorous method for guaranteeing that the perturbations carried out on the system satisfy these constraints. 

Notably, there does exist a fairly established literature on methods that suppose the existence of a parametric model approximation $g_{m,j} ({\bf u}, {\boldsymbol \theta}) \approx g_j ({\bf u})$, define $\Theta$ as the uncertainty set to which ${\boldsymbol \theta}$ belongs, and then attack the problem via probabilistic formulations by ensuring that $g_{m,j} ({\bf u},{\boldsymbol \theta}) \leq 0$ with sufficiently high probability \citep{Kall1994,Zhang2002,Sahinidis2004,Li2008,Quelhas:12}. However, while such methods are theoretically just and robust, they suffer from four major practical drawbacks:

\begin{enumerate}[(i)]
\item the requirement of a parametric model,
\item the restriction that the uncertainty be parametric, and that $\Theta$ be known,
\item the computational issues that arise with probabilistic constraints,
\item the conservatism that results from the probabilistic constraints reducing the set of admissible ${\bf u}$.
\end{enumerate}

\noindent Drawback (i) becomes debilitating when the system at hand is difficult to model, while (ii) is more problematic since many employed models are, often by practical requirement, simplifications and thereby prone to \emph{structural} errors \citep{Chachuat2009}. Drawback (iii) is likely to be significant when the models have many decision variables, many uncertain parameters, and are involved. Simplifications, such as linearizing the model with respect to ${\boldsymbol \theta}$ \citep{Zhang2002}, may be used to avoid this, but ultimately come with the loss of rigor that one would expect from an approximation. Finally, (iv) can be extremely problematic when the parametric uncertainty set is large -- as may often occur in practice \citep{Li2008,Quelhas:12} -- since this may limit the perturbation options, with only a small collection of ${\bf u}$ being deemed ``safe''.

The methodology proposed in the present work avoids these difficulties while maintaining the rigor. Taking a model-free, back-off approach, we simplify and generalize the results of \cite{Bunin:SCFOImp} to derive positive values, $b_{j}$, that, for a given ${\bf u}^*$, allow us to state the guarantee

\begin{equation}\label{eq:imply}
g_{j} ({\bf u}^*) \leq -b_{j} \Rightarrow g_{j} ({\bf u}) \leq 0, \;\; \forall {\bf u} \in \mathcal{B}_e,
\end{equation}

\noindent where 

$$\mathcal{B}_e = \{ {\bf u} : \| {\bf u} - {\bf u}^*  \|_2 \leq \delta_e \}.$$ 

\noindent Verbally, this means that given a decision-variable set ${\bf u}^*$ known to satisfy the constraints with some slack, one is able to provide a guarantee that the entire ball of radius $\delta_e$ surrounding ${\bf u}^*$ will satisfy the constraints as well, thereby allowing the user to perturb anywhere within this ball without fear of constraint violation. Despite being local, such a result is nevertheless very useful as it allows a high degree of freedom -- a ball permitting perturbation sets of any geometry. As will be shown, the value $b_{j}$ will depend on the local sensitivities of $g_{j}$ around ${\bf u}^*$, but can nevertheless be computed without requiring much effort from the user. Conversely, $\delta_e$ is the \emph{sole tuning parameter} set by the user and represents, in some sense, the magnitude of perturbation considered as ``sufficiently exciting'' for identification given the particular problem.

To date, this result has already been integrated into the SCFO experimental optimization solver \citep{SCFOug}, where it is used to ensure accurate linear and quadratic regression, but it is expected that the generality of the result make it applicable to many algorithms and contexts. In this paper, its usefulness is illustrated for a much simpler optimization algorithm -- the evolutionary operation (EVOP) method of Box \citep{Box1957, Box:69}. As the original method searches to maximize an experimental function by perturbing in a hypercube around the best known ${\bf u}$, it is made coherent with the result here by ensuring that the cube lie inside $\mathcal{B}_e$, with ${\bf u}^*$ then defined as the best known reference point. By forcing ${\bf u}^*$ to always satisfy (\ref{eq:imply}), it thus follows that all exploration by the modified EVOP version satisfy the constraints.

The remainder of this paper is organized as follows. The required mathematical concepts and the derivation of the appropriate constraint back-offs are presented in Section 2. Section 3 then provides a robust extension of (\ref{eq:imply}) that accounts for noise/error in the function values, together with a general discussion of potential implementation issues. The constrained EVOP algorithm is presented in Section 4, and its effectiveness is illustrated for three case-study problems. Section 5 concludes the paper.  

\section{Derivation of Sufficient Back-Offs}
\label{sec:deriv}

So as to keep the forthcoming analysis relatively simple, the following assumption on the continuity and differentiability of $g_j$ is made.

\begin{ass}\label{assum:C2}
The functions $g_{j}$ are continuously differentiable ($\mathcal{C}^1$) on an open set containing $\mathcal{B}_e$.
\end{ass}

This then allows for the definition of bounds on the sensitivities of $g_j$.

\begin{define}\label{def:lip}
The local Lipschitz constants of $g_{j}$ are defined as any constants $\kappa_{ji}$ satisfying

\begin{equation}\label{eq:lipcon}
- \kappa_{ji} \leq \frac{\partial g_{j}}{\partial u_i} \Big |_{{\bf u}} \leq \kappa_{ji}, \; \forall {\bf u} \in \mathcal{B}_e.
\end{equation}

\end{define}

The existence of these constants follows from Assumption \ref{assum:C2} and the boundedness of $\mathcal{B}_e$. They may be used to bound the violation of a given $g_j$ via the local \emph{Lipschitz upper bound}.

\begin{lemma}\label{lem:lipbound}
Let ${\bf u}_a, {\bf u}_b \in \mathcal{B}_e$. It follows that

\begin{equation}\label{eq:lipbound}
g_{j} ({\bf u}_b) \leq g_{j} ({\bf u}_a) + \displaystyle \sum_{i=1}^{n_u} \kappa_{ji} |u_{b,i} - u_{a,i}|.
\end{equation}

\end{lemma}
\begin{pf}
See \cite{Bunin:Lip}. \qed
\end{pf}

Finally, the Lipschitz bound may be exploited by substituting ${\bf u}_a \rightarrow {\bf u}^*$ and ${\bf u}_b \rightarrow {\bf u}$ in (\ref{eq:lipbound}) to generate a \emph{Lipschitz polytope} around the point ${\bf u}^*$.

\begin{define}\label{def:lippoly}
Let $\mathcal{L}_{j}$ denote the Lipschitz polytope of the constraint $g_{j}$ around ${\bf u}^*$, defined as the set

\begin{equation}\label{eq:lippoly}
\mathcal{L}_{j} = \left\{ {\bf u} : g_{j} ({\bf u}^*) + \displaystyle \sum_{i=1}^{n_u} \kappa_{ji} | u_{i} - u_{i}^*| \leq 0 \right\}.
\end{equation}

\end{define}

The Lipschitz polytope has two important properties that should be apparent by inspection:

\begin{enumerate}[(i)]
\item ${\bf u} \in \mathcal{L}_{j} \cap \mathcal{B}_e \Rightarrow g_{j} ({\bf u}) \leq 0$, 
\item $g_{j}({\bf u}^*) \leq 0 \Rightarrow \mathcal{L}_{j} \neq \varnothing$.
\end{enumerate}

\noindent The ``double membership'' of ${\bf u} \in \mathcal{L}_j$ and ${\bf u} \in \mathcal{B}_e$ in (i) is required to ensure both that ${\bf u}$ satisfies the upper bound of (\ref{eq:lipbound}) and that the bound itself is valid to begin with, respectively. Property (ii) should be evident if one just considers ${\bf u} := {\bf u}^*$ when $g_j ({\bf u}^*) \leq 0$.

Furthermore, it is clear that the content (hypervolume) of $\mathcal{L}_{j}$ increases monotonically as $g_{j} ({\bf u}^*)$ decreases -- i.e., $\mathcal{L}_{j}$ admits more and more implementable points because of the terms $| u_i - u_i^* |$ being allowed to grow larger while satisfying the inequality.

\begin{figure}
\begin{center}
\includegraphics[width=10cm]{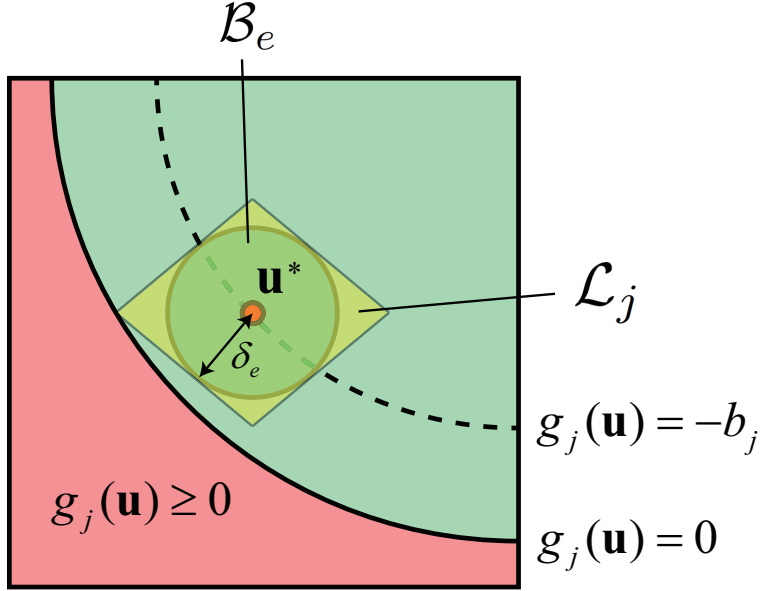}
\caption{Geometric illustration of the concept behind using back-offs to enforce the existence of an excitation ball $\mathcal{B}_e$. The back-off is chosen large enough so that a ball of radius $\delta_e$ may be \emph{inscribed} inside the Lipschitz polytope corresponding to the back-off.}
\label{fig:lipball}
\end{center}
\end{figure}

It is this observation that inspires the foundations of the present work, illustrated geometrically in Figure \ref{fig:lipball}. If $g_{j} ({\bf u}^*)$ can be forced to remain sufficiently low -- i.e., if $g_{j} ({\bf u}^*)$ can be made to satisfy the constraint with a certain back-off -- then one can always guarantee the existence of a non-empty Lipschitz polytope centered at ${\bf u}^*$. Additionally, because the ball $\mathcal{B}_e$ is also centered at ${\bf u}^*$, it may be inscribed inside the polytope, with Property (i) above then sufficient to guarantee that all points in this ball satisfy the constraint $g_{j} ({\bf u}) \leq 0$.

The required size of this back-off is easily derived by exploiting the Cauchy-Schwarz inequality.

\begin{thm}\label{thm:main} Let ${\boldsymbol \kappa}_{j}$ denote the vector $(\kappa_{j1},...,\kappa_{jn_u})$ corresponding to a given $\delta_e$ and ${\bf u}^*$. If $b_{j} \geq \delta_e \| {\boldsymbol \kappa}_{j} \|_2$, then the implication (\ref{eq:imply}) holds.
\end{thm}

\begin{pf} Defining the quantity $\Delta u_i = | u_i - u_i^* |$ and the vector $\Delta {\bf u} = \left( \Delta u_1, \hdots, \Delta u_{n_u} \right)$, let us restate the definition of the Lipschitz polytope as

$$
\mathcal{L}_{j} = \left\{ {\bf u} : g_{j} ({\bf u}^*) + {\boldsymbol \kappa}_{j}^T \Delta {\bf u} \leq 0 \right\}.
$$

From the Cauchy-Schwarz inequality, it follows that 

\begin{equation}\label{eq:cauchy}
{\boldsymbol \kappa}_{j}^T \Delta {\bf u} \leq \| {\boldsymbol \kappa}_j \|_2 \| \Delta {\bf u} \|_2 = \| {\boldsymbol \kappa}_j \|_2 \| {\bf u} - {\bf u}^* \|_2,
\end{equation}

\noindent which then allows us to consider a (ball) subset of ${\mathcal{L}_j}$:

$$
\mathcal{\hat L}_{j} = \left\{ {\bf u} : g_{j} ({\bf u}^*) + \| {\boldsymbol \kappa}_j \|_2 \| {\bf u} - {\bf u}^* \|_2 \leq 0 \right\},
$$

\noindent where $\mathcal{\hat L}_j \subseteq \mathcal{L}_j$ is evident since any ${\bf u}$ satisfying

$$
g_j({\bf u}^*) + \| {\boldsymbol \kappa}_j \|_2 \| {\bf u} - {\bf u}^* \|_2 \leq 0
$$

\noindent must also satisfy

$$
g_j({\bf u}^*) + {\boldsymbol \kappa}_j^T \Delta {\bf u} \leq 0
$$

\noindent by virtue of (\ref{eq:cauchy}).

It is then sufficient to show that $\mathcal{B}_e \subseteq \mathcal{\hat L}_j$ for the back-off specified. Because every ${\bf u} \in \mathcal{B}_e$ satisfies $\| {\bf u} - {\bf u}^* \|_2 \leq \delta_e$, it follows that any point in $\mathcal{B}_e$ will belong to $\mathcal{\hat L}_j$ if the inequality $g_j ({\bf u}^*) + \delta_e \| {\boldsymbol \kappa}_j \|_2 \leq 0$ holds, which is precisely what is ensured by the back-off. Thus, $\mathcal{B}_e \subseteq \mathcal {\hat L}_j \subseteq \mathcal {L}_j$, with Property (i) of the Lipschitz polytope then leading to the desired result. \qed

\end{pf}

It is natural to ask what one should do if, given $\delta_e$, ${\bf u}^*$, and (consequently) $\kappa_{ji}$, the condition $g_j ({\bf u}^*) \leq -\delta_e \| {\boldsymbol \kappa}_j \|_2$ does not hold. One solution consists in choosing a different, safer ${\bf u}^*$ from the available measurements, maintaining the same $\delta_e$, and hoping that the back-off is satisfied for this new choice. Alternatively, one may keep ${\bf u}^*$ the same and reduce $\delta_e$. Assuming that $g_j ({\bf u}^*) < 0$, it is not hard to show that for $\delta_e$ sufficiently small the back-off will be met. Consider, for example, the reduction sequence 

$$
\left\{ \frac{\delta_e}{2^n} \right\}, \; n \in \mathbb{N} = \left\{ 0,1,2,\hdots \right\},
$$

\noindent together with the corresponding $\mathcal{B}_e$:

$$
\left\{ {\bf u} : \| {\bf u} - {\bf u}^*  \|_2 \leq \frac{\delta_e}{2^n} \right\}.
$$

\noindent Clearly, each reduced $\mathcal{B}_e$ will be a strict subset of the previous, from which we may deduce that the $\kappa_{ji}$ valid for $n := 0$ will remain valid for all $n$, thus allowing us to upper bound the $\| {\boldsymbol \kappa}_j \|_2$ portion of the back-off and to conclude that $\delta_e \| {\boldsymbol \kappa}_j \|_2$ may be made arbitrarily close to 0 by choosing $n$ sufficiently large ($\delta_e$ sufficiently small), thereby resulting in a back-off that is satisfied by $g_j ({\bf u}^*)$.

\section{Implementation Issues}
\label{sec:implement}

On the surface, the derived result seems easy to implement as it simply states that if, given ${\bf u}^*$, one wants to perturb anywhere in $\mathcal{B}_e$ without fear of constraint violation, one only needs to confirm that $g_{j} ({\bf u}^*) \leq -\delta_e \| {\boldsymbol \kappa}_j \|_2, \; \forall j$. If this condition does not hold, one can then choose a different ${\bf u}^*$ from the set of applied ${\bf u}$ or decrease $\delta_e$. However, there remain a number of implementation aspects that merit discussion.

\subsection{Setting of Lipschitz Constants}

Theorem \ref{thm:main} provides a robust guarantee only if it is assumed that the Lipschitz constants provided satisfy (\ref{eq:lipcon}). While one can argue that picking very large, conservative values is sufficient, this has the obvious performance drawback of increasing $\delta_e \| {\boldsymbol \kappa}_j \|_2$ and thus $b_j$, which could make satisfying $g_{j} ({\bf u}^*) \leq -b_j$ impossible. As such, one would prefer setting these constants in more intelligent ways.

Apart from some scarce and limited results in the global optimization literature \citep{Hansen1992, Strongin1973, Wood1996}, the problem of estimating Lipschitz constants effectively is still very open. Some techniques have been proposed in the recently submitted manuscript of \cite{Bunin2016Lip}, and  suggest mixing \emph{a priori} knowledge about the functions $g_j$, which may provide good initial guesses of the constants, together with data-driven refinements, making explicit use of the fact that for a sufficiently local or approximately linear region the Lipschitz constants differ little from the local function derivatives. In the case that parametric models $g_{m,j} ({\bf u},{\boldsymbol \theta})$ are available, one could bolster the \emph{a priori} knowledge with the robust model-based estimates

$$
\kappa_{ji} := \mathop {\rm max} \limits_{{\footnotesize \begin{array}{c} {\bf u} \in \mathcal{B}_e \\ {\boldsymbol \theta} \in \Theta \end{array}}} \Bigg |  \frac{\partial g_{m,j}}{\partial u_i} \Big |_{{\bf u},{\boldsymbol \theta}}  \Bigg |
$$

\noindent as the Lipschitz-constant estimates. This mix of techniques has been employed in the SCFO solver \citep{SCFOug} and has generally led to relatively robust performance in test problems, although constraint violations may still occur from time to time \citep{ExpOpt}.

It is worth noting that the EVOP method proposed in Section \ref{sec:evopt} appears to handle this issue \emph{very well}, using iterative local linear regression to estimate the Lipschitz constants in a manner that requires no input from the user and yet avoids constraint violation completely.

\subsection{Accounting for Noise/Error}

Enforcing $g_{j} ({\bf u}^*) \leq -b_j$ also relies on having accurate knowledge of $g_{j} ({\bf u}^*)$. While it is reasonable to expect that in most applications one will be able to either measure or estimate the values of these functions at ${\bf u}^*$, for experimental functions it is generally the case that a perfect measurement or estimation will be impossible, and that there will only be access to the corrupted values, $\hat g_j^*$. A typical assumption \citep{Hotelling1941, Myers2009,Marchetti2010,More2011} is that this corruption is additive -- i.e., that

$$
\hat g_j^* = g_j ({\bf u}^*) + w,
$$

\noindent with $w$ a stochastic element for which the estimates of at least the mean and variance are available. In this case, it becomes possible to compute a high-probability bounding value, $g_j ({\bf u}^*) \leq \overline g_j^*$, using, in the most general case, Chebyshev's inequality \citep{More2011}, or something less conservative if better assumptions on the nature of the noise are available. More involved techniques for computing $\overline g_j^*$ are outlined in Section 4 of \cite{Bunin:SCFOImp} and in the recent work of \cite{Bunin2016Lip}.

One may then work with the robust condition $\overline g_{j}^*\leq-b_j$ instead, since satisfaction of this condition implies the satisfaction of $g_{j} ({\bf u}^*) \leq -b_j$ with a high probability.

\subsection{Accomodating Numerical Constraints}
\label{sec:numcon}

So far, $g_j$ has been treated from a very general perspective, and has been assumed only to be a $\mathcal{C}^1$ function over $\mathcal{B}_e$. However, there are plenty of constraints for which \emph{much more} knowledge is available. In particular, when the constraint $g_{j}$ is a numerical function that can be evaluated by hand or by computer for any desired ${\bf u} \in \mathcal{B}_e$, the guarantee of perturbing in $\mathcal{B}_e$ without incurring violation of the given constraint is fairly easy. Namely, it suffices to ensure that

$$
\mathop {\max} \limits_{{\bf u} \in \mathcal{B}_e} g_j ({\bf u}) \leq 0,
$$ 

\noindent since this trivially guarantees that all perturbations in $\mathcal{B}_e$ are feasible for that constraint, without requiring Lipschitz constants or robust upper bounds. For some cases, it may be that $g_j$ is an involved \emph{nonconcave} function, potentially making its maximization over a ball a numerically challenging (global optimization) task. When this occurs, one could attempt to take a concave upper bound or, if worse comes to worst, estimate its Lipschitz constants -- a significantly easier task for a numerical function -- and then employ the general (albeit conservative) result of Theorem \ref{thm:main}.

A special subclass of numerical constraints that is even easier to work with is that of the bound constraints, $u_i^L \leq u_i \leq u_i^U$, which are extremely relevant since they tend to occur in virtually any well-defined experimental investigation. By simple intuition and quick inspection, one should be able to see that employing a back-off of $\delta_e$ is sufficient for these bounds. Not surprisingly, it is possible to arrive at this same result by Theorem \ref{thm:main}, as the ${\boldsymbol \kappa}_j$ vector for these constraints may simply be taken as the vector of $n_u-1$ zeros with unity in the $i^{\rm th}$ spot, thus leading to $\delta_e \| {\boldsymbol \kappa}_j \|_2 = \delta_e$.

One could, however, choose not to back off from the bound constraints since, given their orthogonal nature, one is always guaranteed to retain a full orthant in which one can perturb. Given with the feasibility of $\mathcal{B}_e$, this means that one could ignore these back-offs and still retain $1/2^{n_u}$ of $\mathcal{B}_e$ for safe perturbation. For most cases, this fraction of $\mathcal{B}_e$ still offers enough room for the user to accomplish what they aim to, such as estimating a derivative or locally exploring the function's behavior.

\subsection{Scaling and the Tuning of the Excitation Radius}

It is evident that the results provided are not scale-invariant, and that a ball of radius $\delta_e$ could offer very poor perturbation in directions that vary over a much greater domain than others -- geometrically, this may be seen as the problem of the efficiency of inscribing a ball inside a Lipschitz polytope. One solution could be to generalize the notion of the ball to an ellipse, and to rederive the result for this more versatile case. However, a simple approach found to work well has been to simply scale the independent variables so that each $u_i$ varies between 0 and 1, in which case perturbing in a ball is reasonable.

This aside, one may still have to make a nontrivial decision about what $\delta_e$ to set. Ultimately, this choice should depend on the particular problem that the user is trying to solve, with all available \emph{a priori} knowledge exploited to yield an appropriate choice. For example, if the goal is to explore as much of the experimental space as possible and to build a model that encompasses a large, perhaps even global, domain, it may be of interest to set $\delta_e$ to be as large as possible. If the goal is to perturb just enough to estimate a derivative while avoiding significant corruption due to noise, one may want to choose $\delta_e$ as a function of the expected signal-to-noise ratio \citep[\S 3.14]{SCFOug}. When information about the function's curvature is available, it may be possible to choose $\delta_e$ so as to strike a balance between the corruption due to noise and the corruption due to nonlinearity \citep{Marchetti2010}. Finally, the \emph{ad hoc}, post-scaling choices of $\delta_e := 0.01, 0.05, 0.10$ have also often worked well, in the author's experience.

\section{Feasible-Side Evolutionary Operation}
\label{sec:evopt}

An evolutionary operation (EVOP) algorithm is constructed for the feasible-side solution of experimental optimization problems having the form

\begin{equation}\label{eq:optprob}
\begin{array}{rll}
\mathop {{\rm{minimize}}}\limits_{\bf{u}} & \phi ({\bf{u}}) & \\
{\rm{subject}}\hspace{1mm}{\rm{to}} & g_{j}({\bf{u}}) \leq 0, & j = 1,...,n_{g} \vspace{1mm}  \\
& u_i^L \leq u_i \leq u_i^U, & i = 1,...,n_u,
\end{array}
\end{equation}

\noindent where both $\phi$ and $g_j$ will, in general, represent experimental functions. As in many experimental optimization problems of this type, there is a constant trade-off between adapting the decision variables ${\bf u}$ so as to minimize $\phi$ and perturbing them so as to learn more about the local or general natures of $\phi$ and $g_j$. Both tasks are necessary -- one cannot optimize without perturbation, but one cannot optimize if one spends all of the experimental resources perturbing the system for knowledge, either.

Two main reasons motivate this choice of application. First, the derived result of Theorem \ref{thm:main} is innately coherent with the direct-search algorithmic nature of EVOP, which iteratively perturbs in a local region and then shifts this region so that it is centered around the best found point. It is thus quite natural to let EVOP perturb in $\mathcal{B}_e$ and then shift ${\bf u}^*$ accordingly, which immediately provides the traditional EVOP with constraint satisfaction guarantees. The second reason has to do with the practical usefulness of the constructed scheme. Ever since its inception in 1957 \citep{Box1957}, EVOP has enjoyed great popularity in industry, largely because of its extreme simplicity. In addition to being simple, it would later be shown to be theoretically well founded as well, with the slight addition of a step-size rule allowing it to enjoy the global convergence properties of Torczon's generalized pattern search algorithms \citep{Torczon1997}. However, both its industrial and theoretical success has, for the most part, been limited to simple problems with bound constraints, and while one could use the standard (i.e., penalty-function) methods \citep[Ch. 14]{Conn2000} \citep[Ch. 12]{Fletcher1987} to convert problems with $g_j$ constraints into the bound-constrained form for which EVOP is readily applicable \citep{Rutten2015}, these approaches are ultimately not safe and can at most only offer the guarantee that the constraints are satisfied upon convergence.

From this point of view, the algorithm proposed in this section is believed to carry great potential as a stand-alone contribution that, in addition to being effective for the problems tested, is also extremely easy to apply and requires minimal input from the user.

\subsection{Description of the Algorithm}

Prior to stating the algorithm, let us first make some knowledge assumptions so as to ensure its basic functionality.

First, it is assumed that, for a given tested $\overline {\bf u}$, one measures both $\phi$ and $g_j$ with additive white Gaussian noise:

$$
\begin{array}{ll}
\hat \phi = \phi (\overline {\bf u}) + w_{\phi}, & w_\phi \sim \mathcal{N} (0,\sigma_\phi^2) \vspace{1mm}  \\
\hat g_j = g_j (\overline {\bf u}) + w_{j}, & w_j \sim \mathcal{N} (0,\sigma_j^2).
\end{array}
$$

\noindent For simplicity, it will be assumed that $\sigma_\phi$ and ${\boldsymbol \sigma} = (\sigma_1,...,\sigma_{n_g})$ are known.

It will also be assumed that a sufficiently ``safe'' initial ${\bf u}^*$ is available, so that the $2n_u$ test points generated by perturbing the individual elements of ${\bf u}^*$ by $\pm \delta_e$,

$$(u_1^* \pm \delta_e,..., u_{n_u}^*),$$
$$(u_1^*, u_2^* \pm \delta_e,..., u_{n_u}^*),$$ 
$$\vdots$$
$$(u_1^*, u_2^*,..., u_{n_u}^* \pm \delta_e),$$

\noindent are, together with ${\bf u}^*$, safe and satisfy the constraints.

Scaled variables will be denoted with $(\tilde \cdot)$ and defined as

$$
\tilde u_i = \frac{u_i - u_i^L}{u_i^U - u_i^L}, \;\;\; i = 1,...,n_u.
$$

\noindent This scaling step will be implicit in the algorithm statement that follows -- i.e., we will switch between the scaled and unscaled variables as needed without explicitly including the affine operation above in the algorithm steps.

\;
\;
\noindent {\bf Algorithm 1 (Feasible-Side EVOP)}
\;

\begin{enumerate}
\item (Initialization) The initial reference point ${\bf u}^*$, the standard deviations $\sigma_\phi$ and ${\boldsymbol \sigma}$, and the decision-variable bounds ${\bf u}^L$, ${\bf u}^U$ are provided. The excitation radius $0.5 \geq \delta_e > 0$ is set by the user.
\item (Evaluation at Reference Point) Apply ${\bf u}^*$ to the experimental system to obtain $\hat \phi^*$ and $\hat g_j^*$. Define $\tilde {\bf U} := (\tilde {\bf u}^*)^T$, ${\boldsymbol \phi} := \hat \phi^*$, ${\bf G} := [\hat g_1^*\; \cdots \; \hat g_{n_g}^*]$. 
\item (Perturbation) For $i = 1,...,n_u$:
\begin{enumerate}
\item (Perturbing by $\pm \delta_e$) Obtain $\tilde {\bf u}_+$ and $\tilde {\bf u}_-$ by perturbing the $i^{\rm th}$ element of $\tilde {\bf u}^*$ by $\delta_e$ and $-\delta_e$, respectively. 
\item (Upper Bound Constraint Check) If $u_{+,i} \leq u_i^U$, apply ${\bf u}_+$ to the experimental system to obtain the corresponding measurements $\hat \phi$ and $\hat g_j$, denoted by $\hat \phi^+$ and $\hat g_j^+$, and augment $\tilde {\bf U}$, ${\boldsymbol \phi}$, and ${\bf G}$:

$$
\tilde {\bf U} := \left[ \begin{array}{c} \tilde {\bf U} \\ {\tilde {\bf u}_+^T}  \end{array} \right], \;\; {\boldsymbol \phi} := \left[ \begin{array}{c} {\boldsymbol \phi} \\ \hat \phi^+ \end{array} \right],
$$

$$
{\bf G} := \left[ \begin{array}{c} {\bf G} \\ \hat g_1^+\; \cdots \; \hat g_{n_g}^+  \end{array} \right].
$$

\item (Lower Bound Constraint Check) If $u_{-,i} \geq u_i^L$, apply ${\bf u}_-$ to the experimental system to obtain the corresponding measurements $\hat \phi$ and $\hat g_j$, denoted by $\hat \phi^-$ and $\hat g_j^-$, and augment $\tilde {\bf U}$, ${\boldsymbol \phi}$, and ${\bf G}$:

$$
\tilde {\bf U} := \left[ \begin{array}{c} \tilde {\bf U} \\ {\tilde {\bf u}_-^T}  \end{array} \right], \;\; {\boldsymbol \phi} := \left[ \begin{array}{c} {\boldsymbol \phi} \\ \hat \phi^- \end{array} \right],
$$

$$
{\bf G} := \left[ \begin{array}{c} {\bf G} \\ \hat g_1^-\; \cdots \; \hat g_{n_g}^-  \end{array} \right].
$$

\item (Set Number of Points Tested) If 

$$0 \leq \tilde u_i^* \pm \delta_e \leq 1,$$

\noindent set $s_i := 2$. Otherwise, set $s_i := 1$.

\end{enumerate}

\item (Local Linear Regression) Fit a linear model, 

$$\beta_0 + \sum_{i=1}^{n_u} \beta_i \tilde u_i,$$ 

\noindent to the observed cost function values, with the $\beta_1,...,\beta_{n_u}$ coefficients recovered as the first $n_u$ elements of the least-squares solution $[\tilde {\bf U}\; {\bf 1}]^\dagger {\boldsymbol \phi}$. Set $\nabla \hat \phi^* = (\beta_1, ..., \beta_{n_u})$ as the estimate of the gradient of $\phi$ at $\tilde {\bf u}^*$. Repeat the same procedure with the columns of ${\bf G}$ to obtain the gradient estimates $\nabla \hat g_j^*$ for $j = 1,...,n_g$.

\item (Lipschitz Constants) For each $i = 1,...,n_u$ and $j = 1,...,n_g$, set 

$$\kappa_{ji} := \Bigg | \frac{\partial \hat g_j^*}{\partial \tilde u_i} \Bigg | + 6 \frac{\sigma_j \sqrt{2}}{s_i \delta_e},$$

\noindent where $\partial \hat g_j^* / \partial \tilde u_i$ denotes the $i^{\rm th}$ element of $\nabla \hat g_j^*$.

\item (Nearly Active Constraints) Define the index set

$$j_A := \left\{ \begin{array}{l} j : \hat g_j + 3\sigma_j \geq -\delta_e \| {\boldsymbol \kappa}_j \|_2 \\ {\rm for\;at\;least\;one\;measurement}\;\hat g_j\;{\rm in\;the}\;j^{\rm th}\;{\rm column\;of}\;{\bf G} \end{array}  \right\} $$
  
\noindent as the index set of nearly active constraints.

\item (Approximation of Lagrangian) Defining the gradient of the Lagrangian at $\tilde {\bf u}^*$ as

$$\nabla L ({\boldsymbol \lambda}) = \nabla \hat \phi^* + \sum_{j=1}^{n_g} \lambda_{j} \nabla \hat g_j^*, $$

\noindent approximate the corresponding Lagrange multipliers, ${\boldsymbol \lambda}^*$, by the solution to the constrained least-squares problem

$$
\begin{array}{rl}
{\boldsymbol \lambda}^* := {\rm arg} \mathop {{\rm{minimize}}}\limits_{\boldsymbol \lambda} & \nabla L ({\boldsymbol \lambda})^T \nabla L ({\boldsymbol \lambda}) \\
{\rm{subject}}\hspace{1mm}{\rm{to}} & \lambda_{j} \geq 0, \;\; \forall j \vspace{1mm}  \\
& \lambda_j = 0, \;\; \forall j \notin j_A.
\end{array}
$$ 

\item (Choose New Reference) Set as $\tilde {\bf u}^*$ the point in $\tilde {\bf U}$ that 

\begin{enumerate}[(i)]
\item has the smallest $\nabla L({\boldsymbol \lambda}^*)^T \tilde {\bf u}$ value,
\item has the corresponding upper bounds $\overline g_j^* := \hat g_{j}^* + 3\sigma_j$ that satisfy $\overline g_j^* \leq -\delta_e \| {\boldsymbol \kappa}_j \|_2, \; \forall j$. 
\end{enumerate}

\noindent If no such $\tilde {\bf u}^*$ exists, maintain the same reference point as before. Define $\tilde {\bf U} := (\tilde {\bf u}^*)^T$, ${\boldsymbol \phi} := \hat \phi^*$, ${\bf G} := [\hat g_1^*\; \cdots \; \hat g_{n_g}^*]$ and return to Step 3.

\end{enumerate}

A number of remarks are in order:

\begin{itemize}
\item As the scaling reduces the variable space to a unit hypercube, there is no need for $\delta_e$ to be set above 0.5, as this would preclude the existence of a feasible $\mathcal{B}_e$.
\item The perturbation scheme used here differs from that of the traditional EVOP methods, which use a $2^{n_u}$ factorial scheme to define the test points \citep{Box1957, Box:69}, and is in this sense more similar to a coordinate search \citep[Ch. 7]{Conn2009}. The rationale for doing this is that requiring $2n_u$ perturbations leads to better efficiency as $n_u$ increases (linear as opposed to geometrical), requires fewer perturbations for $n_u > 2$, and appears to be sufficient for the given tasks.
\item As suggested in Section \ref{sec:numcon}, no back-offs are used for the bound constraints, ensuring instead that they are never violated during the perturbation phase.
\item It is not difficult to show that the derivative estimates obtained by linear regression of the local data set are, as a result of the orthogonality and symmetry of the perturbations, identical to the difference quotients

\vspace{-2mm}
$$
\frac{\partial \hat g_j^*}{\partial \tilde u_i} = \frac{\hat g_j^+ - \hat g_j^-}{2 \delta_e}
$$

\noindent when $s_i: = 2$, or to either

\vspace{-2mm}
$$
\frac{\hat g_j^+ - \hat g_j^*}{\delta_e}\; {\rm or} \; \frac{\hat g_j^* - \hat g_j^-}{\delta_e}
$$

\noindent when $s_i := 1$. It is well known (see, e.g., \cite{Daniel1950}) that the standard deviation for such estimates is equal to $\sigma_j \sqrt{2}/(s_i \delta_e)$. In adding six times this quantity to the absolute value of the derivative estimates, we use half of that amount (``three sigma'') to account for potential errors in the estimates due to noise, and the other half to compensate for the fact that the Lipschitz constants only match the derivatives very locally, and may actually be larger for the $\mathcal{B}_e$ at the future cycle (i.e., the Lipschitz-constants estimation is always one cycle behind). This is important since the new reference is always chosen subject to the restriction that the corresponding $\overline g_j^* \leq -\delta_e \| {\boldsymbol \kappa}_j \|_2, \; \forall j$, where the back-off is defined using the Lipschitz constants at the \emph{current cycle}. This then encourages the guarantee of safe excitation around the new reference.
\item Adding $3\sigma_j$ to the constraint measurements in Steps 6 and 8 robustifies the scheme against measurement noise, resulting in the upper bounding values $\overline g_j := \hat g_j + 3\sigma_j$ that bound the true function values with a high probability of about 99.85\%.
\item Because the problem addressed is constrained and uses Lipschitz bounds to guarantee robust constraint satisfaction, it follows that the scheme could converge suboptimally if allowed to get too close to a constraint. This issue was first observed by \cite{Bunin2011} and dealt with rigorously by \cite{Bunin2013SIAM}, and represents a main drawback of using Lipschitz constants for constraint satisfaction. So as to avoid this issue here, it has been proposed to minimize an approximation of the Lagrangian instead of simply minimizing $\phi$. Note that when ${\bf u}^*$ is far from any constraints, the scheme simply sets ${\boldsymbol \lambda}^* := {\bf 0}$ and the two objectives are equal. However, when some $g_j$ start to get close to activity, the scheme defines the objective as a trade-off between minimizing $\phi$ and lowering the nearly active constraint values, which tends to allow the algorithm to ``slide off'' the nearly active constraints in application. The theoretical rigor of this scheme is not entirely clear, but it may be seen to be consistent upon convergence, in that finding ${\boldsymbol \lambda}^*$ such that $\nabla L({\boldsymbol \lambda}^*) \approx {\bf 0}$ implies that ${\bf u}^*$ is a constrained stationary point, provided that the gradient estimates are not too erroneous.

\end{itemize}

\subsection{Application to Case-Study Examples}

Three problems are chosen from the ExpOpt database of test problems \citep[Problems P2, P3, P6]{ExpOpt} so as to illustrate the performance of the algorithm when applied to realistic case-study scenarios. The first problem consists in varying the feed rate and temperature of a Williams-Otto plant so as to maximize its steady-state profit while honoring an upper limit on a noxious product, and is originally taken from the paper of \cite{Marchetti2013}. The second, adapted from the works of \cite{Gentric1999} and \cite{Francois2005}, deals with minimizing the operating time of a polysterene batch reactor while honoring a minimal molecular weight specification. Here, two ``switching times'' that help define the temperature profile of the batch are taken as the decision variables. Finally, the third problem comes from \cite{Francois2013} and seeks to maximize the steady-state production of a continuous reactor by varying two feed rates subject to two experimental constraints. Two-dimensional problems have intentionally been chosen as they allow for easy visualization and interpretation of the results. The different problem specifications are provided in Table \ref{tab:specs}.

\begin{table}\begin{center}\caption{Problem specifications for the case-study examples.}
\label{tab:specs} 

  \begin{tabular}{|c|c|c|c|c|c|} \hline
    \#  & Initial ${\bf u}^*$ & $\sigma_\phi$ & ${\boldsymbol \sigma}$ & ${\bf u}^L$ & ${\bf u}^U$  \\ \hline
    P2  & $(3.5, 72)$ & $0.5$  & $5 \cdot 10^{-4}$ & $(3, 70)$ & $(6, 100)$  \\ \hline
    P3  & $(242, 945)$ & $60$  & $10^{4}$ & $(50, 600)$ & $(450, 1000)$ \\ \hline
    P6  & $(14.5, 14.9)$ & $0.1$  & $(0.03, 0.03)$ & $(1, 1)$ & $(50, 50)$ \\ \hline
  \end{tabular}

\end{center}
\end{table}

\begin{figure}
\begin{center}
\includegraphics[width=11cm]{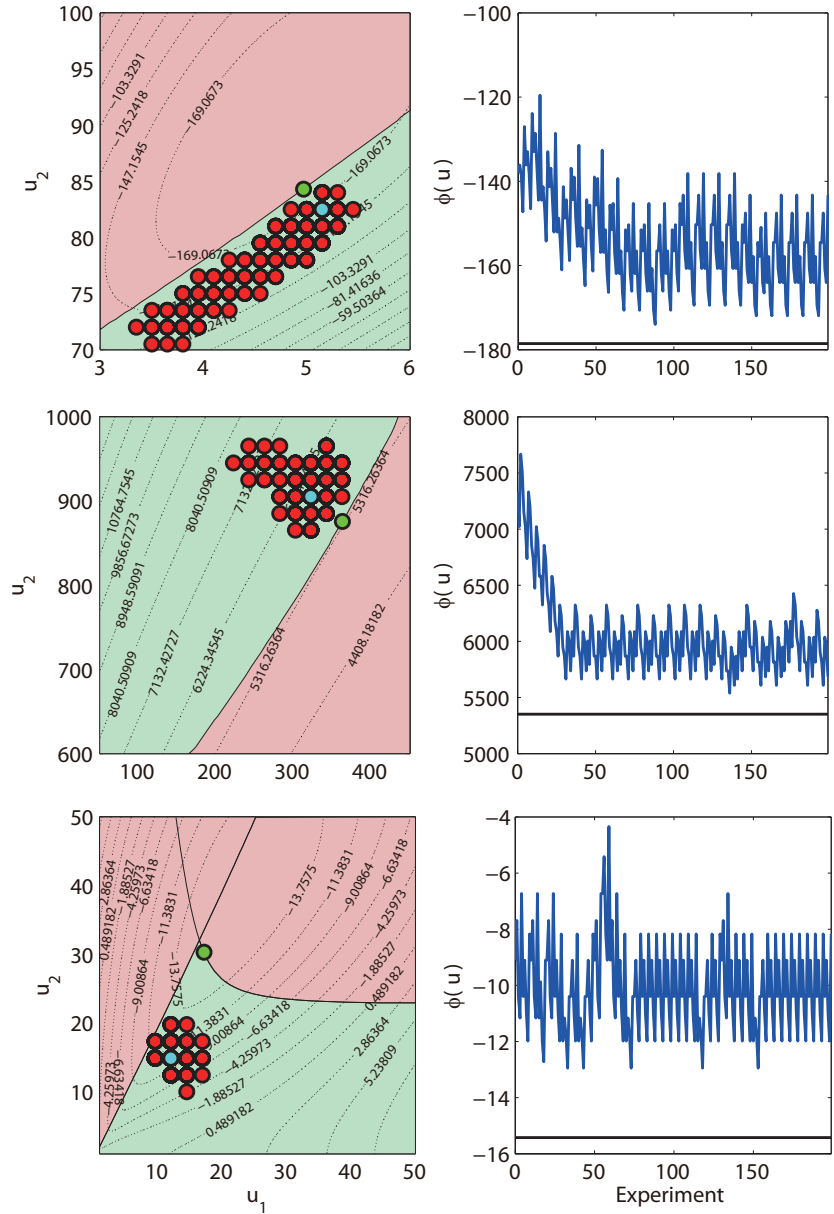}
\caption{Plots of the decision variables and cost function values for Problems P2 (top), P3 (middle), and P6 (bottom). Green points in the decision-variable plots denote the true optimum, red points denote the individual experiments, and blue points denote the reference (best) point at the final tested EVOP cycle. Red regions denote the infeasible (not safe) portions of the operating space as defined by the constraints. In plotting the cost function values, the blue lines denote the true (noiseless) values, while the constant black line denotes the value at the optimum.}
\label{fig:res}
\end{center}
\end{figure}

\begin{figure}
\begin{center}
\includegraphics[width=11cm]{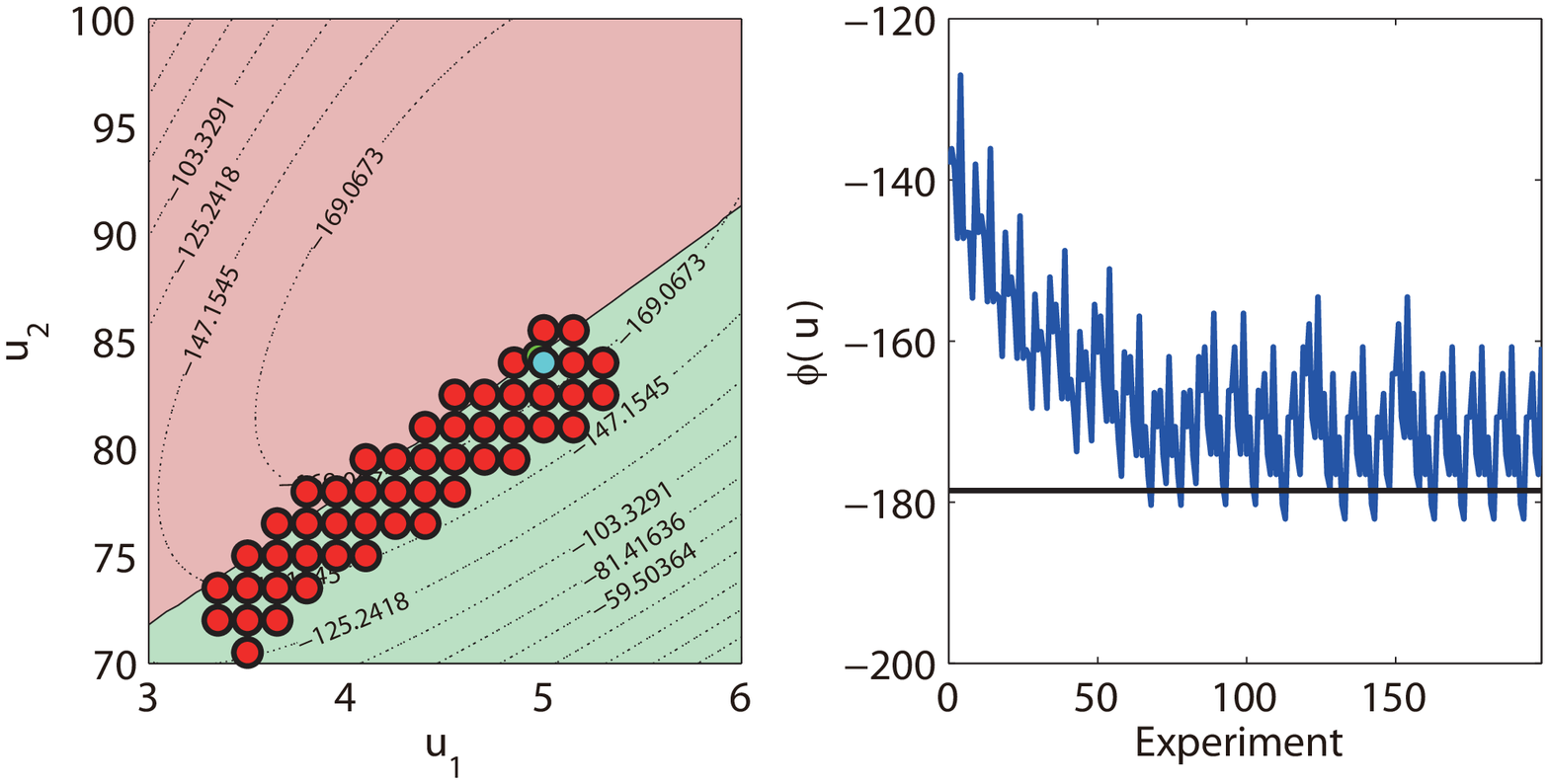}
\caption{Illustration of EVOP performance for Problem P2 when the back-off is not employed.}
\label{fig:resNB}
\end{center}
\end{figure}

Inputting the specifications into Algorithm 1, the algorithm's performance is tested for the \emph{ad hoc} choice of $\delta_e := 0.05$. The results are given in Figure \ref{fig:res} and show that the algorithm succeeds in exploring the decision-variable space while satisfying the constraints -- not a single violation is observed. By contrast, suppose that one tried to apply the same algorithm to Problem P2 but without the application of the safety back-off derived in Theorem \ref{thm:main}. This is done by substituting $-\delta_e \| {\boldsymbol \kappa}_j \|_2 \rightarrow 0$ in Steps 6 and 8 of Algorithm 1. The result is given in
Figure \ref{fig:resNB}. Again, it is seen that the use of the Lagrangian leads to the algorithm successfully converging to a neighborhood of the optimum. However, constraint violations occur repeatedly along the convergence trajectory.

\begin{figure}
\begin{center}
\includegraphics[width=11cm]{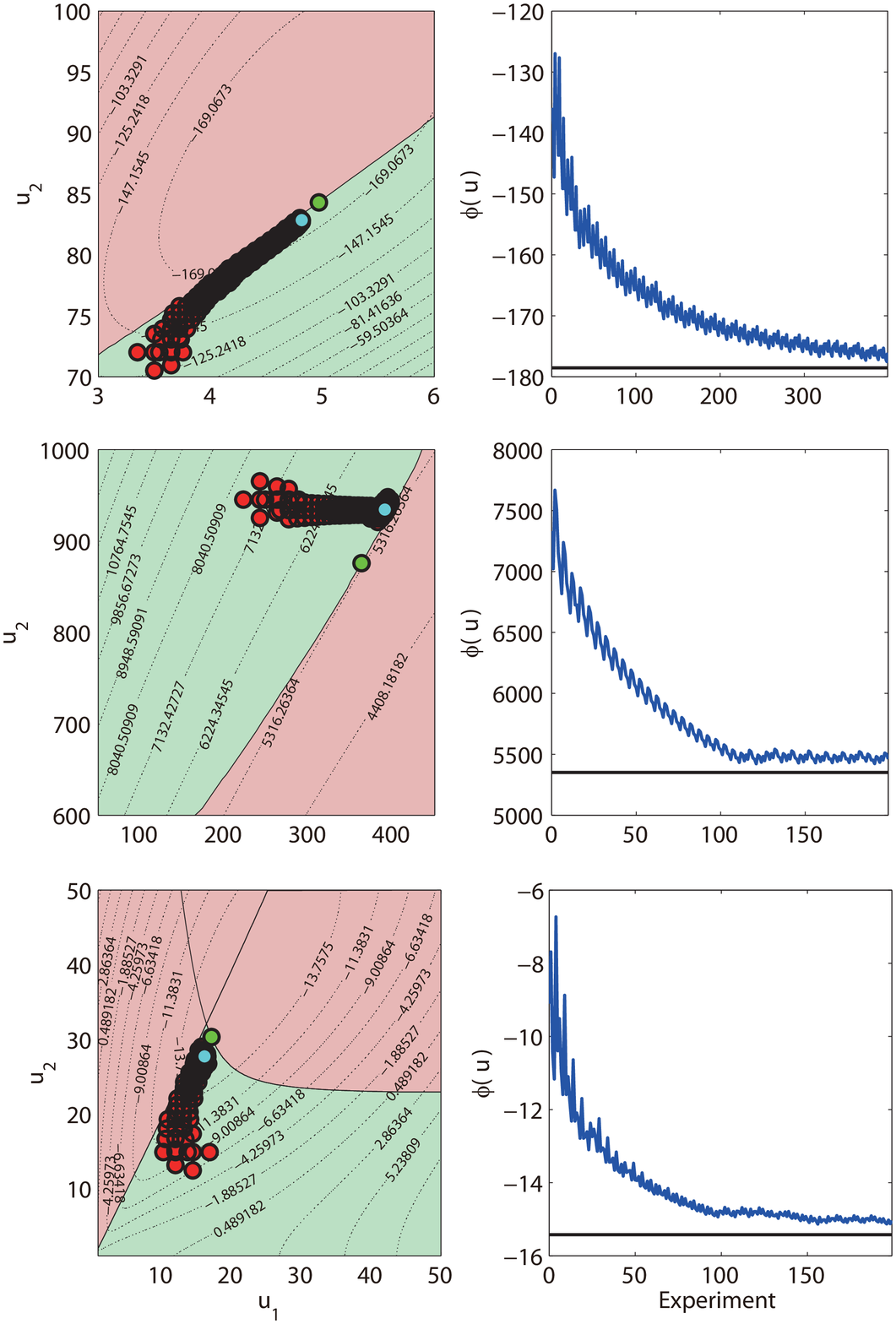}
\caption{Plots of the decision variables and cost function values for Problems P2 (top), P3 (middle), and P6 (bottom) when $\delta_e$, $\sigma_\phi$, and ${\boldsymbol \sigma}$ are gradually reduced over the course of operation.}
\label{fig:resR}
\end{center}
\end{figure}

The convergence properties of the algorithm are of interest -- in particular, one notices that the algorithm converges to a region that is relatively suboptimal in Problem P6. This is due to the geometry of the problem, the conservatism of the back-offs, and the size of the noise in the constraint measurements. In other words, the algorithm cannot make further progress while robustly (to $3\sigma_j$) satisfying the specified back-offs. While proving convergence to an optimum in a mathematically rigorous manner is outside the scope of this paper, for $\mathcal{C}^1$ functions one should expect the algorithm to converge properly if both the noise and back-offs go to 0 asymptotically, as $\delta_e \downarrow 0$ would lead to the linear regression converging to the true function gradient in the absence of noise \citep[\S 2.4]{Conn2009}, with the removal of noise also removing the conservatism introduced by using the upper bounds $\hat g_j+3\sigma_j$. To test this hypothesis empirically, let us use a modified implementation of Algorithm 1 where 

$$\begin{array}{l}
\delta_{e} (k) := \delta_e/\sqrt{k},\\ 
\sigma_{\phi} (k) := \sigma_{\phi}/\sqrt{k}, \\
{\boldsymbol \sigma} (k) := {\boldsymbol \sigma}/\sqrt{k},
\end{array}$$ 

\noindent with $k$ starting at 1 and being incremented by 1 after each cycle of the algorithm. The results are presented in Figure \ref{fig:resR} and largely confirm one's expectations, with feasible-side convergence very close to the optimum achieved for all problems. Although the neighborhood of convergence appears to be far from the optimum in Problem P3, this is due to the geometry of the problem -- because both the cost and constraint functions are close to linear in the region of the optimum, there is a wide range of points that achieve a low cost and are close to stationarity all along the constraint.

\section{Closing Remarks}

This paper has contributed to the problem of obtaining information about an experimental function in the presence of general $\mathcal{C}^1$ constraints, and we have derived a rigorous constraint back-off that, when satisfied at some ${\bf u}^*$, guarantees that the user may perturb anywhere in the ball of radius $\delta_e$ around ${\bf u}^*$ without incurring constraint violations. This result is believed to constitute a useful contribution to scientific problems dealing with constrained experimental spaces, as it offers a straightforward way to guarantee safety while exploring the experimental space.

While the results simplify greatly for numerical constraints, in the case of general (likely experimental) constraints one is forced to obtain local sensitivity bounds -- Lipschitz constants -- for the constraint function in order to compute the appropriate back-off. To do this well may be challenging, but relevant methods do exist and may work quite well in certain contexts \citep{Bunin2016Lip}.

The feasible-side EVOP optimization algorithm has provided an interesting application of the derived back-off result, and has been shown to work very well with respect to constraint satisfaction and convergence for three different case studies. While more involved than the traditional EVOP procedure, the algorithm nevertheless retains its ease of implementation, requiring the user to only set $\delta_e$ prior to applying it to a problem. It goes without saying that numerous performance improvements are possible but have not been the focus of this paper. One could, for example, use \emph{all} of the measurements obtained since initialization to filter out the noise, an idea routinely employed in the SCFO solver \citep[\S 3]{SCFOug} and recently proposed by \cite{Gao2016} in the context of multilayer real-time optimization. One could also attempt to choose search directions in a more intelligent manner, such as what is done in some derivative-free methods \citep{Conn2009,Lewis2000}. Clearly, one could generalize the method to handle noise that is not white Gaussian, as well, although the computational effort may increase due to a lack of closed-form expressions.


\end{document}